\documentclass{gtart}
\def\ifplaintex{\expandafter\ifx\csname documentclass\endcsname\relax}

\ifplaintex 
\else
\fi

\expandafter\ifx\csname beginpicture\endcsname\relax
\expandafter\ifx\csname documentclass\endcsname\relax
\input pictex \else
\input prepictex \input pictex \input postpictex \fi\fi

\def\gt{{\mathsurround=0pt\it $\cal G\mskip-2mu$eometry \&\ 
$\cal T\!\!$opology}}        %  journal title in recommended style

\def\gtp{{\mathsurround=0pt\it $\cal G\mskip-2mu$eometry \&\ 
$\cal T\!\!$opology $\cal P\!$ublications}}  % GT publications

\def\lognumber#1{\def\thelognumber{#1}}
\def\volumenumber#1{\def\thevolumenumber{#1}}
\def\papernumber#1{\def\thepapernumber{#1}}
\def\volumeyear#1{\def\thevolumeyear{#1}}

\def\pagenumbers#1#2{\def\startpage{#1}\def\finishpage{#2}}
\def\published#1{\def\publishdate{#1}}
\def\proposed#1{\def\theproposer{#1}}
\def\seconded#1{\def\theseconders{#1}}
\def\received#1{\def\receiveddate{#1}}
\def\revised#1{\def\reviseddate{#1}}
\def\accepted#1{\def\accepteddate{#1}}

\def\coverauthors#1{\def\thecoverauthors{#1}}
\def\asciiauthors#1{\def\theasciiauthors{#1}}
\def\asciiaddress#1{\def\theasciiaddress{#1}}
\def\asciiemail#1{\def\theasciiemail{#1}}
\long\def\asciiabstract#1{\long\def\theasciiabstract{#1}}

\let\\\par\let\thelognumber\relax
\let\thevolumenumber\relax\let\thepapernumber\relax
\let\thevolumeyear\relax\let\thesamplenumber\relax\let\startpage\relax
\let\finishpage\relax\let\publishdate\relax\let\receiveddate\relax
\let\reviseddate\relax\let\accepteddate\relax\let\theasciititle\relax
\let\theasciiauthors\relax\let\theasciiaddress\relax
\let\theasciiabstract\relax
\let\theasciiemail\relax\let\theshortauthors\relax\let\theshorttitle\relax
\let\thecoverauthors\relax

\long\def\maketitlep{   % start of definition of \maketitlep

\count0=\startpage

\gt\hfill      %   Journal title (top left) 
\beginpicture
\setcoordinatesystem units <0.33truein, 0.33truein> point at 2.2 0.9
\setplotsymbol ({$\cal G$})
\plotsymbolspacing=9truept
\circulararc 315 degrees from 0 1 center at 0 0
\setplotsymbol ({$\cal T$})
\circulararc 315 degrees from 1 -1 center at 1 0
\endpicture
\break
{\small\ifx\thesamplenumber\relax % sample?  
Volume \else Sample
\fi\thevolumenumber\ (\thevolumeyear)
\startpage--\finishpage\nl
Published: \publishdate}
\vglue 0.5truein plus 0.4fil minus 0.1truein

{\parskip=0pt\leftskip 0pt plus 1fil\def\\{\par\smallskip}{\ifplaintex\large
\else\Large\fi\bf\thetitle}\par\medskip}   

\vglue 0pt plus 0.1fil 

{\parskip=0pt\leftskip 0pt plus 1fil\def\\{\par}{\sc\theauthors}
\par\medskip}

\vglue 0pt plus 0.1fil 

{\small\parskip=0pt\let\newline\\
{\leftskip 0pt plus 1fil\def\\{\par}{\sl\theaddress}\par}
\expandafter\ifx\theemail\relax    % email address?
\relax\else\vglue 5pt plus 0.02fil minus 2pt\def\\{\stdspace{\rm 
and}\stdspace} 
\cl{Email:\stdspace\tt\theemail}\fi
\ifx\theurl\relax                  % URL given?
\relax\else\vglue 5pt plus 0.02fil minus 2pt\def\\{\stdspace{\rm 
and}\stdspace}
\cl{URL:\stdspace\tt\theurl}\fi\par}

\vglue 7pt plus 0.3fil minus 3pt

{\bf Abstract}
\vglue 5pt plus 0.1fil minus 2pt

\theabstract

\vglue 7pt plus 0.3fil minus 3pt

{\bf AMS Classification numbers}\quad Primary:\quad \theprimaryclass

Secondary:\quad \thesecondaryclass

\vglue 5pt plus 0.3fil minus 2pt

{\bf Keywords}\quad \thekeywords

\vglue 10pt plus 0.5fil minus 5pt

{\small  Proposed: \theproposer\hfill Received: \receiveddate\nl
Seconded: \theseconders\hfill 
\ifx\reviseddate\relax                         % paper revised?
Accepted: \accepteddate                        % no
\else
Revised: \reviseddate                          % yes
\fi}
\eject
}       %  end of definition of \maketitlep

\let\maketitlepage\maketitlep
\let\maketitle\maketitlepage

\font\phead=cmsl9 scaled 950
\font=cmsl9 scaled 1050
\font\pnum=cmbx10 scaled 913
\font\lnum=cmbx10 
\font\pfoot=cmsl9 scaled 950
\font=cmsl9 scaled 1050
\ifplaintex
\headline{\vbox to 0pt{\vskip -4.5mm\line{\small\phead\ifnum
\count0=\startpage ISSN 1364-0380 (on line)
1465-3060 (printed) \hfill {\pnum\folio}\else\ifodd\count0\def\\{ }% 
\ifx\theshorttitle\relax\thetitle\else\theshorttitle\fi\hfill{\pnum\folio}
\else\def\\{ and }{\pnum\folio}\hfill\ifx\theshortauthors\relax\theauthors
\else\theshortauthors\fi\fi\fi}\vss}}
\footline{\vbox to 0pt{\vglue 0mm\line{\small\pfoot\ifnum\count0=\startpage
\copyright\ \gtp\hfill\else
\gt, Volume \thevolumenumber\ (\thevolumeyear)\hfill\fi}\vss
}}
\else
\makeatletter
\def\@oddhead{{\small\lhead\ifnum\count0=\startpage ISSN 1364-0380 (on line)
1465-3060 (printed) \hfill {\lnum\number\count0}\else\ifodd\count0
\def\\{ }\ifx\theshorttitle\relax \thetitle \else\theshorttitle\fi\hfill
{\lnum\number\count0}\else\def\\{ and }{\lnum\number\count0}
\hfill\ifx\theshortauthors\relax 
\theauthors\else\theshortauthors\fi\fi\fi}}\def\@evenhead{\@oddhead}
\def\@oddfoot{\small\lfoot\ifnum\count0=\startpage\copyright\ \gtp\hfill\else
\gt, Volume \thevolumenumber\ (\thevolumeyear)\hfill\fi}
\def\@evenfoot{\@oddfoot}
\makeatother
\fi

\newwrite\gtoutfile
\long\gdef\makeheadfile{  %%% start of definition of \makeheadfile
{\def\\{, }\def\s{ }
\immediate\openout\gtoutfile head.xxx
\immediate\write\gtoutfile{To: math@arxiv.org}
\immediate\write\gtoutfile{Subject: put or rep NNNNN:pppp}
\immediate\write\gtoutfile{--text follows this line--}
\immediate\write\gtoutfile{Proxy-for: \ifx\theasciiauthors\relax
\theauthors\else\theasciiauthors\fi\s<\ifx\theasciiemail\relax\theemail\else\theasciiemail\fi>}
\immediate\write\gtoutfile{\noexpand\\}
\immediate\write\gtoutfile{Authors: \ifx\theasciiauthors\relax
\theauthors\else\theasciiauthors\fi}
{\def\\{ }\immediate\write\gtoutfile{Title: \ifx\theasciititle\relax
\thetitle\else\theasciititle\fi}}
\immediate\write\gtoutfile{Subj-class: GT or SG or MG etc}
\immediate\write\gtoutfile{MSC-class: \theprimaryclass\ifx\thesecondaryclass\relax\else, \thesecondaryclass\fi}
\immediate\write\gtoutfile{Journal-ref: Geom. Topol. \thevolumenumber
(\thevolumeyear) \startpage-\finishpage}
\immediate\write\gtoutfile{Comments: Published by Geometry and Topology at}
\immediate\write\gtoutfile{\s\s http://www.maths.warwick.ac.uk/gt/GTVol\thevolumenumber/paper\thepapernumber.abs.html}
\immediate\write\gtoutfile{\noexpand\\}
\immediate\write\gtoutfile{}
\ifx\theasciiabstract\relax
\immediate\write\gtoutfile{\theabstract}\else
\immediate\write\gtoutfile{\theasciiabstract}\fi
\immediate\write\gtoutfile{}
\immediate\write\gtoutfile{\noexpand\\}
\immediate\write\gtoutfile{}
\immediate\closeout\gtoutfile}}  %%% end of definition of \makeheadfile

\def\maketitlepage{\maketitlep\makeheadfile}
\let\maketitle\maketitlepage

\lognumber{290}
\volumenumber{7}\papernumber{6}\volumeyear{2003}
\pagenumbers{225}{254}
\received{1 November 2002}
\revised{19 March 2003}
\accepted{20 March 2003}
\published{24 March 2003}
\proposed{John Morgan}
\seconded{Yasha Eliashberg, Robion Kirby}

\usepackage{amsmath,amsfonts,amscd,flafter,rlepsf}

\newtheorem{theorem}{Theorem}[section]

\newtheorem{cor}[theorem]{Corollary}

\newtheorem{lemma}[theorem]{Lemma}

\theoremstyle{remark}

\newtheorem{defn}[theorem]{Definition}

\newcommand\Rk{\mathrm{rk}}

\newcommand{\Q}{\mathbb{Q}}
\newcommand{\R}{\mathbb{R}}

\newcommand{\Z}{\mathbb{Z}}

\newcommand{\OneHalf}{\frac{1}{2}}

\newcommand{\Zmod}[1]{\Z/{#1}\Z}

\newcommand{\cm}{\cdot}

\newcommand{\Nbd}[1]{{\mathrm{nd}}(#1)}
\newcommand{\nbd}[1]{\Nbd{#1}}
\newcommand{\CDisk}{D}

\newcommand{\ModSWfour}{\mathcal{M}}
\newcommand{\ModFlow}{\ModSWfour}

\newcommand{\SpinC}{{\mathrm{Spin}}^c}

\newcommand\abuts\Rightarrow
\newcommand\Sym{\mathrm{Sym}}

\newcommand\spinccanf{k}

\newcommand\relspinc{\underline{\spinc}}

\newcommand\Filt{\mathcal F}

\newcommand\x{\mathbf x}

\newcommand\y{\mathbf y}

\newcommand\ModSphere{\ModFlow\left({\mathbb S}\longrightarrow 
\Sym^{g-1}(\Sigma_{1})\times \Sym^2(\Sigma_{2})\right)}
\newcommand\ModSpheres\ModSphere
\newcommand\CF{CF}

\newcommand\CFa{\widehat{CF}}
\newcommand\CFp{\CFb}
\newcommand\CFm{\CF^-}

\newcommand\HFp{\HFb}

\newcommand\CFinf{CF^\infty}
\newcommand\HFinf{HF^\infty}
\newcommand\CFb{CF^+}
\newcommand\HFa{\widehat{HF}}
\newcommand\HFb{HF^+}
\newcommand\gr{\mathrm{gr}}
\newcommand\Mas{\mu}
\newcommand\UnparModSp{\widehat \ModSp}
\newcommand\UnparModFlow\UnparModSp
\newcommand\Mod\ModSp

\newcommand{\cald}{{\mathcal D}}

\newcommand{\spinc}{\mathfrak s}

\newcommand\Real{\mathrm Re}

\newcommand\ModMaps{\mathcal M}
\newcommand\ModSp\ModMaps

\newcommand\Ta{{\mathbb T}_{\alpha}}
\newcommand\Tb{{\mathbb T}_{\beta}}

\newcommand\Dual{\mathcal D}
\newcommand\Duality\Dual

\newcommand\InjMod[1]{{\mathcal T}^+_{#1}}
\newcommand\Vertices{\mathrm{Vert}}

\newcommand\CFK{CFK}
\newcommand\HFK{HFK}

\newcommand\CFKa{\widehat\CFK}

\newcommand\CFKinf{\CFK^{\infty}}

\newcommand\HFKa{\widehat\HFK}

\newcommand\Mark{m}

\newcommand\FiltPt{z}
\newcommand\BasePt{w}

\newcommand\States{\mathfrak S}

\begin{document}

\title{Heegaard Floer homology and alternating knots}

\author{Peter Ozsv\'ath\\Zolt{\'a}n Szab{\'o}}
\coverauthors{Peter Ozsv\noexpand\'ath\\Zolt{\noexpand\'a}n Szab{\noexpand\'o}}
\asciiauthors{Peter Ozsvath, Zoltan Szabo} 
\address{Department of
Mathematics, Columbia University\\New York 10027, USA}
\email{petero@math.columbia.edu}
\secondaddress{Department of
Mathematics, Princeton University\\New Jersey 08540, USA}
\secondemail{szabo@math.princeton.edu}

\asciiaddress{Department of
Mathematics, Columbia University\\New York 10027, USA\\and\\Department of
Mathematics, Princeton University\\New Jersey 08540, USA}
\asciiemail{petero@math.columbia.edu, szabo@math.princeton.edu}

\begin{abstract}  
In an earlier paper, we introduced a knot invariant for a
null-homologous knot $K$ in an oriented three-manifold $Y$, which is
closely related to the Heegaard Floer homology of $Y$.  In this paper
we investigate some properties of these knot homology groups for knots
in the three-sphere. We give a combinatorial description for the
generators of the chain complex and their gradings.  With the help of
this description, we determine the knot homology for alternating
knots, showing that in this special case, it depends only on the
signature and the Alexander polynomial of the knot (generalizing a
result of Rasmussen for two-bridge knots). Applications include new
restrictions on the Alexander polynomial of alternating knots.
\end{abstract} 

\asciiabstract{In an earlier paper, we introduced a knot invariant for a
null-homologous knot K in an oriented three-manifold Y, which is
closely related to the Heegaard Floer homology of Y.  In this paper
we investigate some properties of these knot homology groups for knots
in the three-sphere. We give a combinatorial description for the
generators of the chain complex and their gradings.  With the help of
this description, we determine the knot homology for alternating
knots, showing that in this special case, it depends only on the
signature and the Alexander polynomial of the knot (generalizing a
result of Rasmussen for two-bridge knots). Applications include new
restrictions on the Alexander polynomial of alternating knots.}


\maketitlepage
\section{Introduction} 

In~\cite{HolDisk}, we introduced a collection of Abelian groups
associated to a closed, oriented three-manifold $Y$, the Heegaard
Floer homology of $Y$.\footnote{In most of this paper, we work with
the simplest version of Heegaard Floer homology, $\HFa(Y)$, for
three-manifolds with $H_1(Y;\Z)=0$.  In this case $\HFa(Y)$ is a
$\Z$-graded, finitely generated Abelian group.}  In~\cite{Knots}, we
introduced a ``knot filtration'' on the Heegaard Floer homology of a
three-manifold $Y$ which is induced from a null-homologous knot $K$ in
$Y$, see also~\cite{Rasmussen}. Taking the homology groups of the
associated graded object, we obtain Floer homology groups
$\HFKa(Y,K)=\bigoplus_{i\in\Z}\HFKa(Y,K,i)$ (where the integer $i$
appearing here corresponds to the filtration level) which are
topological invariants of the knot $K$.

Our aim here is to study these invariants in the case where the
ambient three-manifold is the three-sphere, in which case
$\HFa(S^3)=\Z$.  Working with a suitable Heegaard diagram compatible
with a planar projection of the knot, we describe the ``classical''
aspects of the Floer theory -- generators of the knot complex, their
filtration levels, and absolute gradings -- in terms of combinatorics
of a (generic) knot projection to the plane (though the differentials
in the knot complex still elude such a description).  With this
combinatorial description in hand, we are able to completely determine
the Heegaard Floer homology for alternating knots -- i.e.\  those which
admit a projection for which the crossing types alternate between
overcrossings and undercrossings -- and give some topological
applications. More calculations based on these descriptions will be
given in a future paper~\cite{calcKT}.

\subsection{Classical Floer data for classical knots}

Let $K\subset S^3$ be a knot. To define the knot Floer homology,
we must work with a Heegaard diagram for $S^3$ which is compatible with
the knot $K$. Specifically, we require that the knot $K$ is supported
entirely inside one of the two handlebodies, meeting 
exactly one of the attaching disks in a single transverse intersection point.
As in~\cite{Knots}, we then obtain a set of generators $X$ 
of the chain complex $\CFa(S^3)$ (the chain complex whose homology calculates
$\HFa(S^3)\cong \Z$), which are endowed with a pair of
integer-valued functions, the filtration level $\Filt$ 
and the absolute grading $\gr$.

For the boundary operator on $\CFa(S^3)$, if $\x$ corresponds to a
generator and $\y$ appears with non-zero multiplicity in the expansion
of $\partial \x$, then $\gr(\y)=\gr(\x)-1$, while $\Filt(\x)\geq
\Filt(\y)$.  Thus, the associated graded complex for the 
filtration $\Filt$, which we denote here by $\CFKa(S^3,K)$, is also
freely generated by $X$, but its boundary operator now preserves $\Filt$,
and hence we have the splitting
$$\CFKa(S^3,K)=\bigoplus_{\{i\in\Z\}}\CFKa(S^3,K,i),$$ where
$\CFKa(S^3,K,i)$ is the subcomplex generated by elements $\x\in X$
with $\Filt(\x)=i$.
We reiterate: although this associated chain complex $\CFKa(S^3,K)$
depends on the choice of Heegaard diagram used for $S^3$, its homology
$$\HFKa(S^3,K)=\bigoplus_{i\in\Z}\HFKa(S^3,K,i)$$ does not.

A {\em decorated projection} for $K$ is a generic knot projection of
$K$, together with a choice of a distinguished edge.
In Section~\ref{sec:Heegaard}, we associate a natural Heegaard diagram
for $K$ to any decorated knot projection for $K$. This allows us to
describe the generators $X$ and the functions $\Filt$ and $\gr$ in
terms of the knot projection.  For the description of the generators,
we use the notion of {\em states} introduced by Kauffman for the
Alexander polynomial, see~\cite{Kauffman}. We recall this briefly
here.

Let $K\subset S^3$ be an oriented knot, and fix a decorated projection
of $K$. The projection gives a planar graph $G$ whose vertices
correspond to the double-points of the projection of $K$. Since $G$ is
four-valent, there are four distinct quadrants (bounded by edges)
emanating from each vertex, each of which is a corner of the closure
of some region of $S^2-G$.  Let $m$ denote the number of vertices of
$G$. Clearly, $G$ divides $S^2$ into $m+2$ regions. The two regions
which share the distinguished edge  will be denoted $A$ and
$B$. In fact, we will always choose our projections so that $A$ is the
unbounded region.

\begin{defn}
A {\em Kauffman state} (cf.~\cite{Kauffman}) for a decorated knot
projection of $K$ is a map which associates to each vertex of $G$ one
of the four in-coming quadrants, so that:
\begin{itemize}
\item the quadrants associated to distinct vertices are subsets of distinct
regions in $S^2-G$
\item none of the quadrants is a corner of the distinguished regions
$A$ or $B$.
\end{itemize}
\end{defn}

It is easy to see that a Kauffman state sets up a one-to-one
correspondence between vertices of $G$ and the connected components of
$S^2-G-A-B$. There is a very simple description of Kauffman states in
graph-theoretic terms (see also~\cite{Kauffman}).  The regions in the
complement of the planar projection can be colored black and white in
a chessboard pattern, by the rule that any two regions which share an
edge have opposite color. There is then an associated ``black graph'',
whose vertices correspond to the regions colored black, and whose
edges correspond to vertices in $G$ which connect the opposite black
regions. In these terms, Kauffman states are in one-to-one
correspondence with the maximal subtrees of the black graph (under a
correspondence which associates to a Kauffman state $x$ the union of
vertices of $G$, thought of now as edges in the black graph, to which
$x$ associates a black quadrant).

Let ${\mathfrak S}$ denote the set of Kauffman states for our
decorated knot projection. We define two functions
$\Filt\colon\States\longrightarrow \Z$ and $M\colon
\States\longrightarrow \Z$, called the filtration level and absolute
grading respectively.

\begin{figure}[ht!]
\cl{\epsfxsize3in\epsfbox{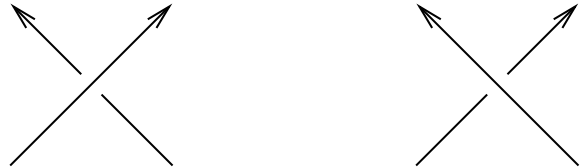}}
\caption{\label{fig:CrossingSigns}
{\bf{Crossing conventions}}\qua Crossings of the first kind are assigned
$+1$, and those of the second kind are assigned $-1$.}
\end{figure}

\begin{figure}[ht!]
\cl{\epsfxsize3in\epsfbox{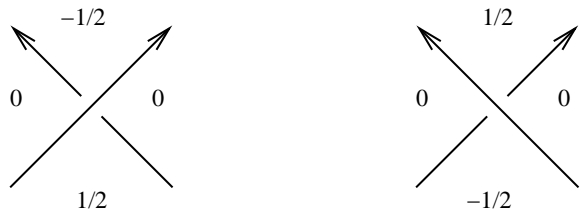}}
\caption{\label{fig:PointTowards}
{\bf{Local filtration level contributions $s(x,v)$}}\qua We have
illustrated the local contributions of $s(x,v)$ for both kinds of
crossings. (In both pictures, ``upwards'' region is the one which the
two edges point towards.)}
\end{figure}

To describe the filtration level, note that the orientation on the
knot $K$ associates to each vertex $v\in G$ a distinguished quadrant
whose boundary contains both edges which point towards the vertex
$v$. We call this the quadrant which is ``pointed towards'' at $v$.
There is also a diagonally opposite region which is ``pointed away
from'' (i.e.\ its boundary contains the two edges pointing away from
$v$).  We define the local filtration contribution of $x$ at $v$,
denoted by $s(x,v)$, by the following rule (illustrated in
Figure~\ref{fig:PointTowards}), where $\epsilon(v)$ denotes the sign
of the crossing (which we recall in Figure~\ref{fig:CrossingSigns}):
$$2 {\epsilon(v)} s(x,v)= \left\{
\begin{array}{rl}
1 & {\text{$x(v)$ is the quadrant pointed towards at $v$}}\\ -1 &
{\text{$x(v)$ is the quadrant away from at $v$}}
\\ 0 & {\text{otherwise.}}
\end{array}
\right.$$
The filtration level associated to a Kauffman state, then, is given by the sum
$$S(x)=\sum_{v\in\Vertices(G)}s(x,v).$$ Note that the function $S(x)$
is the $T$-power appearing for the contribution of $x$ to the
symmetrized Alexander polynomial, see~\cite{Alexander},
\cite{KauffmanTwo}.

\begin{figure}[ht!]
\cl{\epsfxsize3in\epsfbox{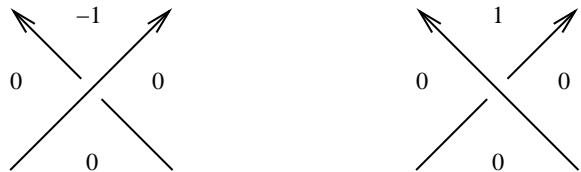}}
\caption{\label{fig:Grading}
{\bf{Local grading contributions $m(x,v)$}}\qua We have
illustrated the local contribution of $m(x,v)$.}
\end{figure}

The grading $M(x)$ is defined analogously.  First, at each vertex $v$,
we define the local grading contribution $m(x,v)$.  This local
contributions is non-zero on only one of the four quadrants -- the
one which is pointed away from at $v$.  At this quadrant, the grading
contribution is minus the sign $\epsilon(v)$ of the crossing, as
illustrated in Figure~\ref{fig:Grading}.  Now, the
grading $M(x)$ of a Kauffman state $x$ is defined by the formula
$$M(x)=\sum_{v\in \Vertices(G)}m(x,v).$$
With these objects in place, we can now state the following:

\begin{theorem}
\label{thm:States}
Let $K$ be a knot in the three-sphere, and choose a decorated knot
projection of $K$.  Then, there is a Heegaard diagram for $K$ with the
property that the knot complex $\CFKa(S^3,K)$ is freely generated by
Kauffman states of the decorated projection.  Moreover, if $\x$
denotes the generator of $\CFKa(S^3,K)$ and $x$ is its corresponding
Kauffman state, then $\Filt(\x)=S(x)$ and $\gr(\x)=M(x)$.
\end{theorem}

\subsection{Heegaard Floer homology for alternating knots}
In the special case where the knot is alternating,
Theorem~\ref{thm:States} easily determines the knot homology
completely in terms of the Alexander polynomial $\Delta_K$ of the knot
$K$ and its signature $\sigma(K)$.

In the following theorem, we use the sign conventions according to
which the signature of the left-handed trefoil is $+2$,
cf.~\cite{Lickorish}.

\begin{theorem}
\label{thm:KnotHomology}
Let $K\subset S^3$ be an alternating knot in the three-sphere, and
write its symmetrized Alexander polynomial as $$\Delta_K(T) =
a_0+\sum_{s>0} a_s (T^s+T^{-s}).$$ Then, $\HFKa(S^3,K,s)$ is supported
entirely in dimension $s+\frac{\sigma}{2}$, and indeed
$$\HFKa(S^3,K,s)\cong \Z^{|a_s|}.$$
\end{theorem}

It is very suggestive to compare the above result with the
corresponding theorem of Lee on the Khovanov homology for alternating
knots (see~\cite{EunSooLee}, see also~\cite{Khovanov},
\cite{BarNatan}, \cite{Garoufalidis}).
We can also use the above theorem to calculate the Heegaard Floer
homologies of three-manifolds obtained by zero-surgeries along the
knot.  In particular, we obtain the following result, which is a
generalization of a theorem of Rasmussen (see~\cite{Rasmussen}) which
calculates the Heegaard Floer homology of three-manifolds obtained as
integer surgeries along two-bridge knots. (In effect, our methods show
that alternating knots are ``perfect'' in Rasmussen's sense.)

To state our results in a useful form, we recall more aspects of the
Heegaard Floer homology package from~\cite{HolDisk}. Specifically,
there is a variant of the Heegaard Floer homology $\HFp(Y)$ which is a
module over the ring $\Z[U]$. This module is related to the variant
$\HFa(Y)$ considered earlier by a canonical long exact sequence $$
\begin{CD}
\cdots@>>>\HFa(Y)@>>>\HFp(Y)@>{U}>>\HFp(Y)@>>> \cdots
\end{CD}
$$ 
The group $\HFp(Y)$ can be given an absolute $\Zmod{2}$ grading. In fact,
these groups obtain some additional structure, depending on the
the homological properties of $Y$, which we describe now in the case
where $H_1(Y;\Z)\cong \Z$ (as is the case for zero-framed surgery on
a knot in $S^3$). In this case, there is a splitting $\HFp(Y)\cong
\bigoplus_{s\in\Z} \HFp(Y,s)$ with the property that $\HFp(Y,s)\cong
\HFp(Y,-s)$.  
Indeed, the summand $\HFp(Y,s)$ is endowed with a relative
$\Zmod{2s}$
grading\footnote{Let $S$ be an Abelian group. An {\em (absolutely) $S$-graded
Abelian group} is an Abelian group $G$ generated by a set $X$, equipped
with a map $\gr\colon X \longrightarrow S$. A {\em relatively $S$-graded
Abelian group} is an Abelian group $G$ generated by a set $X$, equipped
with a relative grading $\gr'\colon X\times X \longrightarrow S$
with 
$\gr'(x,y)+\gr'(y,z)=\gr'(x,z)$.} (compatible
with the $\Zmod{2}$ mentioned earlier), and in particular when $s=0$,
it gives a relative $\Z$-grading.  Indeed the summand with $s=0$ can
be given an absolute grading with values in the set $\OneHalf + \Z$.
Although this choice might seem unnatural at first glance, it fits
neatly with the four-dimensional theory, cf.~\cite{AbsGraded}.
Action by the ring element $U$ decreases all of these gradings by $2$.

Let $\InjMod{k}$ denote the $\Q$-graded
$\Z[U]$-module which is abstractly isomorphic to $\Z[U,U^{-1}]/\Z[U]$,
graded so that multiplication by $U$ decreasing grading by two,
and its bottom-most homogeneous generator has degree $k\in\Q$. Recall that
$\HFp(S^3)\cong \InjMod{0}$.

Given a knot $K$ and an integer $s$, let $t_s(K)$ denote the 
{\em torsion coefficients} defined by
$$t_s(K)=\sum_{j=1}^\infty j a_{|s|+j}$$ (where here the $a_s$
are the coefficients of the symmetrized
Alexander polynomial of $K$). These integers are closely related to
the Milnor torsion of the knot~\cite{MilnorTorsion}, see also~\cite{Turaev}.

Finally, for
$\sigma\in 2\Z$, we let $\delta(\sigma,s)$ be the integer defined by
\begin{equation}
\label{eq:DefOfDelta}
\delta(\sigma,s)=\max(0,\lceil \frac{|\sigma|-2|s|}{4}\rceil).
\end{equation}
Note that $\delta(\sigma,s)$ is the $s^{th}$
torsion coefficient of the $(2,|\sigma|+1)$
torus knot.

\begin{theorem}
\label{thm:FloerHomology}
Let $K$ be an alternating knot, oriented so that
$\sigma=\sigma(K)\leq 0$, and let $S^3_0(K)$ denote the
three-manifold obtained by zero-surgery on $K$.  Then,
\begin{itemize}
\item  for all 
$s>0$, we have  a $\Z[U]$-module isomorphism
$$\HFp(S^3_0(K),s)\cong \Z^{b_s}\oplus
\left(\Z[U]/U^{\delta(\sigma,s)}\right),$$ 
where the first summand is supported in degree
$s+\frac{\sigma}{2}\pmod{2}$, while the second summand has odd parity,
and $\delta(\sigma,s)$ is defined as in Equation~\eqref{eq:DefOfDelta},
\item for $s=0$, we have an isomorphism of graded $\Z[U]$ modules
$$\HFp(S^3_0(K))\cong \Z^{b_0} \oplus \InjMod{-1/2}\oplus 
\InjMod{-2\delta(\sigma,0)+\OneHalf}$$
and the cyclic
summand $\Z^{b_0}$ lies in degree $\frac{\sigma-1}{2}$.
\end{itemize}
Thus, in both cases, $b_s$ is given by the formula
\begin{equation}
\label{eq:GetBs}
(-1)^{s+\frac{\sigma}{2}}b_s=\delta(\sigma,s)-t_s(K).
\end{equation}
\end{theorem}

In fact, Theorem~\ref{thm:FloerHomology} is a formal consequence of
Theorem~\ref{thm:KnotHomology} (for any knot which satisfies the
conclusion of Theorem~\ref{thm:FloerHomology}, the conclusion of
Theorem~\ref{thm:KnotHomology} holds). In fact, there are some
non-alternating knots which satisfy the
conclusion of Theorem~\ref{thm:FloerHomology}. An example is given
in Section~\ref{sec:AltLinks}.

%whose Heegaard diagrams can be further
%simplified, to see that they also satisfy these hypotheses.

%Theorem~\ref{thm:KnotHomology} admits also a straightforward generalization
%to the case of 

\subsection{Applications to the topology of alternating knots}

We now describe some of the consequences of the above calculations for
alternating knots, combined with other results on Heegaard Floer homology.

As a first consequence, we obtain the following calculation of the
``correction terms'' for three-manifolds obtained as surgery on $S^3$
along $K$, for any alternating knot $K$.  This correction term is
defined using the absolute grading on the $\Z[U]$ module $\HFp(Y)$
when $Y$ is an integer homology three-sphere, described
in~\cite{AbsGraded}. (As explained in that reference, when
$H_1(Y;\Z)=0$, the group $\HFp(Y)$ can be endowed with an absolute
$\Z$-grading.)  This number is the analogue of the gauge-theoretic
invariant of Fr{\o}yshov introduced in~\cite{Froyshov},
\cite{FroyshovInst}, constraining the intersection forms of
four-manifolds which bound $Y$.  Specifically, according to
Theorem~\ref{AbsGraded:intro:IntForm} of~\cite{AbsGraded}, if $Y$ is
an integer homology three-sphere, then for each negative-definite
four-manifold $W$ which bounds $Y$, we have the inequality
\begin{equation}
\label{eq:IntFormConstraint}
	\xi^2+\Rk(H^2(W;\Z))\leq 4d(Y),	
\end{equation}
for each characteristic vector $\xi$ for the intersection form
$H^2(W;\Z)$.\footnote{Recall that for a bilinear form $Q$ over a
lattice $V\cong \Z^n$, a {\em characteristic vector} is a vector
$c\in V$ with the property that $Q(v,v)\equiv Q(c,v)\pmod{2}$ for all
$v\in V$.} Recall that Elkies~\cite{Elkies} proves that for any
negative-definite, unimodular form over $\Z$,
$$\max_{\{{\text{characteristic
vectors $\xi$ for $V$}}\}} \xi^2+\Rk(V)\geq 0,$$ with equality holding 
if and
 only if the bilinear form $V$ is diagonalizable (over $\Z$).  In view of
these results, then, $d(Y)$ bounds the ``exoticness'' of intersection
forms of smooth, definite four-manifolds which bound $Y$, providing a
relative version of Donaldson's diagonalization
theorem~\cite{Donaldson}, see also~\cite{Froyshov} and~\cite{FroyshovInst}.

In general, calculating $d(Y)$ is challenging. But for surgeries on
alternating knots, we have the following explicit result
(compare~\cite{Rasmussen}, \cite{FroyshovII}):

\begin{cor}
\label{cor:Froyshov}
Let $K\subset S^3$ be an alternating knot, then
$$d(S^3_1(K))=2\min(0,-\lceil\frac{-\sigma(K)}{4}\rceil).$$
\end{cor}

In another direction, Theorem~\ref{thm:FloerHomology} can be used to give restrictions on
the Alexander polynomials of alternating knots. A classical result of
Crowell and Murasugi (see~\cite{Crowell} and \cite{MurasugiAlt})
states that the coefficients of the symmetrized Alexander polynomial for such a
knot alternate in sign (indeed, the sign of $a_s$ is 
$(-1)^{s+\frac{\sigma}{2}}$). Theorem~\ref{thm:FloerHomology} in turn
immediately gives the following inequality for the torsion
coefficients, which is easily seen to generalize this alternating
phenomenon:

\begin{cor}
\label{cor:AlternateEstimate}
Let $K$ be an alternating knot in the three-sphere. Then for all
$s\in\Z$, we have that
$$(-1)^{s+\frac{\sigma}{2}}(t_s(K)-\delta(\sigma,s))\leq 0,$$ where
$\delta(\sigma,s)$ are the constants defined in
Equation~\eqref{eq:DefOfDelta}.
\end{cor}

For example, consider the nine-crossing knot $K$ appearing in the
standard knot tables as $9_{42}$, see for
example~\cite{BurdeZieschang}. This knot has
\begin{eqnarray*}
\sigma(K)=2 &{\text{and}}& 
\Delta_K(T)=-1+2(T+T^{-1})-(T^2+T^{-2}),
\end{eqnarray*}
i.e.\ its Alexander polynomial is alternating, but it fails
to satisfy the conditions of Corollary~\ref{cor:AlternateEstimate},
so it is not alternating. (Note this particular result
is classical, see~\cite{CrowellTwo}.)

Other restrictions on the Alexander polynomials of alternating knots
have been conjectured by Fox, see~\cite{Fox} (see
also~\cite{Murasugi}, where these properties are verified for a large
class of alternating knots). Specifically, Fox conjectures that for an
alternating knot, the absolute values of the coefficients of the
Alexander polynomial $|a_s|$ are non-increasing in $s$, for $s\geq 0$.
It is easy to see that the above corollary verifies Fox's conjecture
for alternating knots of genus $2$.  For a general alternating knot,
the inequalities stated above for coefficients $s<g-1$ are independent
of Fox's prediction. However, for the first coefficient change,
i.e.\ when $s=g-1$, the above inequalities translate into the following
stronger bound: $$|a_{g-1}|\geq 2|a_{g}|+\left\{
\begin{array}{rl}
-1 & {\text{if $|\sigma|=2g$}} \\
1 & {\text{if $|\sigma|=2g-2$}} \\
0 & {\text{otherwise.}}
\end{array}\right.$$
Finally, we describe a relationship between Theorem~\ref{thm:FloerHomology}
and contact geometry. Recall that a
fibered knot $K\subset S^3$ endows $S^3$ with an open book
decomposition, and hence a contact structure, using a construction of
Thurston and Winkelnkemper, see~\cite{ThurstonWinkelnkemper}. One can
ask which contact structure this is.

For this purpose, recall that a contact structure in $S^3$ has a
classical invariant, the ``Hopf invariant'' of the induced two-plane
field $h(\xi)\in \Z$, which is an integer which uniquely specifies the
homotopy class of $\xi$. This number is defined by $$4
h(\xi)=c_1(\spinccanf)^2+2-2\chi(W)-3\sigma(W),$$ where $W$ is any
almost-complex four-manifold which bounds $S^3$ so that the induced
complex tangencies on its boundary coincide with $\xi$, $\spinccanf$
is the canonical class of the almost-complex structure, $\chi(W)$ is
the Euler characteristic of $W$, and $\sigma(W)$ is the signature of its intersection form.

Using results on the knot homology of fibered knots described
in~\cite{HolDiskContact} (which, in turn, are based on the important
work of Giroux~\cite{Giroux}), we obtain the following:

\begin{cor}
\label{cor:AltContact}
Let $K\subset S^3$ be an alternating, fibered knot of genus $g$, and
let $\xi_K$ denote its induced contact structure over $S^3$. Then,
\begin{equation}
\label{eq:EquationForHopf}
h(\xi_K)=-\frac{\sigma(K)}{2}-g(K).
\end{equation}
Moreover, the induced contact
structure on $S^3$ is  tight if and only
if $h(\xi_K)=0$.
\end{cor}

In~\cite{EliashbergSThree}, Eliashberg classifies contact structures
over $S^3$, showing that for each non-zero integer $i$, there is a
unique contact structure $\xi_i$ whose Hopf invariant is $i$, while
there are two contact structures with vanishing Hopf invariant: the
``standard'' (tight) contact structure, and another 
(overtwisted) one.
Combining Corollary~\ref{cor:AltContact} with Eliashberg's
classification,
we obtain the following:

\begin{cor}
\label{cor:AltContact2}
The standard contact structure and all other contact structures in $S^3$ with 
negative Hopf invariant are precisely those contact structures which
are represented by alternating, fibered knots.
\end{cor}

\subsection{Alternating links}

Theorem~\ref{thm:KnotHomology}, together with many of its
consequences, admits a straightforward generalization to the case of
non-split, alternating links. We state and prove the generalization in
Section~\ref{sec:AltLinks}. As an illustration,
we use this as a stepping-stone for
a calculation of the knot homology for a non-alternating knot, $9_{48}$.

\medskip
\noindent{\bf{Acknowledgements}}\qua The authors wish to warmly thank
Jacob Rasmussen and Andr{\'a}s Stipsicz for interesting conversations.

PSO was supported by NSF grant number DMS 9971950 and a Sloan 
Research Fellowship;
ZSz was supported by NSF grant number DMS 0107792
and a Packard Fellowship.

\section{Proof of Theorem~\ref{thm:States}.}
\label{sec:Heegaard}

We prove here the state-theoretic interpretation of the classical
Floer data, stated in Theorem~\ref{thm:States}. The main ingredient is
a Heegaard diagram which is naturally associated to a decorated knot
projection. (Note that this is not the usual diagram induced from
placing the knot into ``bridge position''.)  We describe this Heegaard
diagram, after briefly recalling some of the ingredients of
the knot Floer complex for knots (specializing for simplicity to the 
case where the ambient manifold is $S^3$).

\subsection{Classical Floer data.}
We give here a rapid description of the data for the knot Floer
complex captured in Theorem~\ref{thm:States}.  We refer the reader
to~\cite{Knots} for a more detailed discussion.

Fix an oriented knot $K\subset S^3$.  A {\em marked Heegaard diagram}
is triple of data
$$(\Sigma,\{\alpha_1,\ldots,\alpha_g\},\{\beta_1,\ldots,\beta_g\},\Mark),$$
where here
\begin{itemize}
\item $\Sigma$ is an oriented surface of genus $g$, 
\item $\{\alpha_1,\ldots,\alpha_g\}$ are pairwise disjoint, embedded
circles in $\Sigma$ representing the attaching circles for a
handlebody $U_\alpha$; similarly,  $\{\beta_1,\ldots,\beta_g\}$ are
pairwise disjoint, embedded circles in $\Sigma$ representing the attaching circles
for a handlebody $U_\beta$
\item $\Mark$ is a marked point on the attaching circle $\beta_1$, which is
disjoint
from the $\alpha_i$
\item the Heegaard diagram describes $S^3$, i.e.\ 
we have a diffeomorphism $U_\alpha\cup_{\Sigma} U_\beta\cong S^3$
\item under this identification, the knot $K$ is supported entirely inside
$U_\beta$, and it is disjoint from the attaching disks for $\beta_j$ with $j>1$, meeting the attaching disk for $\beta_1$ in a single positive, transverse
intersection point.
\end{itemize}
We  consider
the $g$-fold symmetric product $\Sym^g(\Sigma)$, with two
distinguished tori 
\begin{eqnarray*}
\Ta=\alpha_1\times\cdots\times\alpha_g&{\text{and}}&
\Tb=\beta_1\times\cdots\times\beta_g.
\end{eqnarray*} The generators $X$ for the chain
complex $\CFKa(S^3,K)$ are intersection points between $\Ta$ and $\Tb$
in $\Sym^g(\Sigma)$. Let $\BasePt$ and $\FiltPt$ be two points in
$\Sigma$ which are near $\Mark$, but which lie on either side of
$\beta_1$. A choice of orientation on $K$ gives an ordering on these
two points. More precisely, we can find arc $\delta$ connecting
$\FiltPt$ to $\BasePt$ so that $\delta$ is disjoint from all
$\alpha_i$ and $\beta_j$ with $j>1$, meeting $\beta_1$ in a single,
transverse intersection point. Orienting $\delta$ in the same
direction as $K$, we order the two base points so that $\delta$ goes
from $\FiltPt$ to $\BasePt$.

By simple topological considerations
(cf.\ Section~\ref{HolDisk:sec:TopPrelim} of~\cite{HolDisk}), given
two intersection points $\x,\y\in\Ta\cap\Tb$, we can find a Whitney
disk $\phi$ for $\Ta$ and $\Tb$ which connects $\x$ to $\y$; i.e.\  a
map $\phi$ from the standard complex disk $\CDisk$ into
$\Sym^g(\Sigma)$ with the properties that
\begin{eqnarray*}
u\{\zeta \big| \Real(\zeta)\geq 0~{\text{and}}~|\zeta|=1\}\subset \Ta,
&&
u\{\zeta\big| \Real(\zeta)\leq 0~{\text{and}}~|\zeta|=1\}\subset \Tb,\\
u(-i)=\x, &&
u(i)=\y.
\end{eqnarray*}
Let $p$ be a point in
$\Sigma-\alpha_1-\cdots-\alpha_g-\beta_1-\cdots-\beta_g$, we let $n_p(\phi)$
denote the algebraic intersection number of $\phi$ with the
submanifold $p\times \Sym^{g-1}(\Sigma)$. In fact, by choosing
reference points in each connected component of
$\Sigma-\alpha_1-\cdots-\alpha_g-\beta_1-\cdots-\beta_g$, we obtain a
function from this set of regions to $\Z$, denoted $\cald(\phi)$, and
called the {\em domain} associated to $\phi$.  As explained in
Section~\ref{HolDisk:sec:TopPrelim} of~\cite{HolDisk}, when $g>2$, the
homotopy class of $\phi$ is uniquely determined by $\cald(\phi)$. We
denote the set of homotopy classes of Whitney disks by $\pi_2(\x,\y)$.

We now describe the functions $\Filt$ and $\gr$, referred to in the
introduction. The definitions we sketch here are somewhat simpler than
the general definitions given in~\cite{Knots}\footnote{In fact, in~\cite{Knots}, the
filtration is denoted $\relspinc_\Mark$, and in the case we are considering
now, it takes values in the
set of  $\SpinC$ structures on $S^3_0(K)$. This is related to the present
function $\Filt$ by the formula $$\Filt(\x)=\langle
c_1(\relspinc_\Mark(\x)),[{\widehat F}]\rangle,$$ where here ${\widehat F}$ 
is the closed surface in $S^3_0(K)$ obtained by capping off  a Seifert 
surface $F$ for $K$ which respects the orientation of $K$.}, owing to the fact that
our ambient manifold is $S^3$.

First, we discuss $\gr$.  Given $\x$ and $\y$, let $\phi$ be a Whitney
disk connecting $\x$ to $\y$.  We claim that if $\Mas(\phi)$ denotes
the Maslov index of $\phi$, then $\Mas(\phi)-2n_\BasePt(\phi)$ is
independent of the choice of $\phi$. Indeed
(cf.\ Equation~\eqref{HolDisk:eq:RelativeGrading} of~\cite{HolDisk}),
the function $\gr$ is determined up to an additive constant by the
relation $$\gr(\x)-\gr(\y)=\Mas(\phi)-2n_\BasePt(\phi).$$ The
remaining indeterminacy is removed using the Heegaard Floer homology
of $S^3$: $\CFa(S^3)$ is a chain complex which is generated by the
same intersection points $X$ (though its boundary operator allows for
more differentials than the knot Floer complex $\CFKa(S^3,K)$), its
boundary operator decreases $\gr$ by one, and its homology is $\Z$,
supported in a single dimension. The indeterminacy, then, is removed
by the convention that the homology is supported in dimension zero.

Similarly, if $\x$ and $\y$ are a pair of intersection points and
$\phi$ is a Whitney disk connecting them, then we claim that the
difference of intersection numbers $n_\FiltPt(\phi)-n_\BasePt(\phi)$
is independent of the choice of $\phi$.  Indeed
(cf.\ Lemma~\ref{Knots:lemma:RelSpinCThree} of~\cite{Knots}), the
function $\Filt$ is determined up to an additive constant, by the
equation $$\Filt(\x)-\Filt(\y)=n_\FiltPt(\phi)-n_{\BasePt}(\phi).$$
The remaining indeterminacy of $\Filt$ is removed with the help of the
observation that $$\sum_{\x\in\Filt(x)} (-1)^{\gr(\x)} T^{\Filt(\x)}
=T^c \Delta_K(T),$$ where here $\Delta_K(T)$ denotes the symmetrized
Alexander polynomial of $K$ (see Equation~\eqref{Knots:eq:EulerChar}
of~\cite{Knots}), and $c$ is some
integer. The additive indeterminacy of $\Filt$, then, is removed by
requiring that $c=0$.

We have not defined here the differential on the chain
complex. Loosely speaking, the differential on $\CFa(S^3)$ counts
pseudo-holomorphic Whitney disks $u$ with $n_\BasePt(u)=0$, while the
differential on $\CFKa(S^3,K)$ counts those which satisfy
$n_\BasePt(u)=0=n_\FiltPt(u)$. Details are given in~\cite{Knots}.

\subsection{The Heegaard diagram belonging to a decorated knot projection}

Let $K$ be an oriented knot in $S^3$, and fix a decorated knot
projection of $K$ -- i.e.\ a generic planar projection $G$ of $K$ with
$n$ crossings, and a choice of distinguished edge $e$ which appears in
the closure of the unbounded region $A$. With these choices, we
construct a Heegaard diagram for $(S^3, K)$ as follows.

\begin{figure}[ht!]
\relabelbox\small
\cl{\epsfxsize3.2in\epsfbox{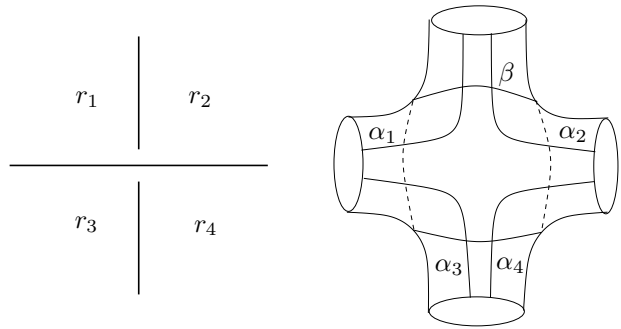}}
\relabel {r1}{$r_2$}
\relabel {r2}{$r_3$}
\relabel {r3}{$r_4$}
\relabel {r4}{$r_1$}
\adjustrelabel <-2pt,0pt> {a1}{$\alpha_3$}
\adjustrelabel <-2pt,0pt> {a2}{$\alpha_4$}
\relabel {a3}{$\alpha_2$}
\relabel {a4}{$\alpha_1$}
\relabel {b}{$\beta$}
\endrelabelbox
\caption{\label{fig:SpecialHeegaard}
{\bf{Special Heegaard diagram for knot crossings}}\qua At each crossing
as pictured on the left, we construct a piece of the Heegaard surface
on the right (which is topologically a four-punctured sphere). The
curve $\beta$ is the one corresponding to the crossing on the left;
the four arcs $\alpha_1,\ldots,\alpha_4$ will close up. (Note that if one
of the four regions $r_1,\ldots,r_4$ contains the distinguished edge $e$,
its corresponding $\alpha$-curve should {\em not} be included).
Note that the Heegaard surface is oriented from the outside.}
\end{figure}

Let $B$ denote the other region which contains the edge $e$, and let
$\Sigma$ be the boundary of a regular neighborhood in $S^3$ of $G$
(i.e.\ it is a surface of genus $n+1$); we orient $\Sigma$ as $\partial
(S^3-\nbd{G})$. We associate to each region $r\in R(G)-A$, an
attaching circle $\alpha_r$ (which follows along the boundary of
$r$). To each crossing $v$ in $G$ we associate an attaching circle
$\beta_v$ as indicated in Figure~\ref{fig:SpecialHeegaard}. In
addition, we let $\mu$ denote the meridian of the knot, chosen to be
supported in a neighborhood of the distinguished edge $e$.

Each vertex $v$ is contained in four (not necessarily distinct)
regions. Indeed, it is clear from Figure~\ref{fig:SpecialHeegaard},
that in a neighborhood of each vertex $v$, there are at most four
intersection points of $\beta_v$ with circles corresponding to these
four quadrants. (There are fewer than four intersection points with
$\beta_v$ if $v$ is a corner of the unbounded region $A$.) Moreover,
the circle corresponding to $\mu$ meets the circle $\alpha_B$ in a
single point (and is disjoint from the other circles).  It is easy to
see that for any choice of marked point $\Mark\in\mu$, the
construction we have just described gives a marked Heegaard diagram
for $K$.

The correspondence between states and generators for the knot complex\break
$\CFKa(S^3,K)$, $\Ta\cap\Tb$, should now be clear: an intersection
point gives at each vertex (i.e.\ $\beta$-curve) one of four quadrants
(corresponding to the up to four $\alpha$-curves). Moreover, since
the meridian $\mu$ meets exactly one $\alpha$-curve, the curve
corresponding to the region $B$, that corresponding $\alpha$-curve is
not assigned to any of the vertices. 

Before turning to the other aspects of Theorem~\ref{thm:States}, we
recall a very useful technical device: the ``Clock Theorem'' of
Kauffman.

\begin{defn}
\label{def:Consecutive}
Two distinct states $x$ and $y$ are said to differ by a {\em
transposition}
if there is a pair of vertices $v_1$ to $v_2$ 
with the property that:
\begin{itemize}
\item $x|_{G-v_1-v_2}\equiv y|_{G-v_1-v_2}$,
\item there is a straight path $P$ from $v_1$ to $v_2$ 
(i.e.\ a path which follows 
along the knot $K$) which does not contain the distinguished edge, so
that $x(v_1)$ and $y(v_1)$ are the two quadrants which contain the
first edge in $P$, and $x(v_2)$ and $y(v_2)$ are the two quadrants
which contain the last edge in $P$.
\end{itemize}
\end{defn}

\begin{theorem}[Kauffman]
Any two distinct states $x$ and $y$ can be connected by a sequence of
transpositions.
\end{theorem}

The proof of the above result can be found in Chapter~2
of~\cite{Kauffman}.

We now establish some lemmas used in the proof of Theorem~\ref{thm:States}.

\begin{lemma}
\label{lemma:SpinCLemma}
Suppose that the two states $x$ and $y$ differ by a transposition, and
let $\x$ and $\y$ denote the corresponding generators of the knot
complex. Then, 
$\Filt(\x)-\Filt(\y)=S(x)-S(y)$.
\end{lemma}

\begin{figure}[ht!]
\cl{\epsfxsize4in\epsfbox{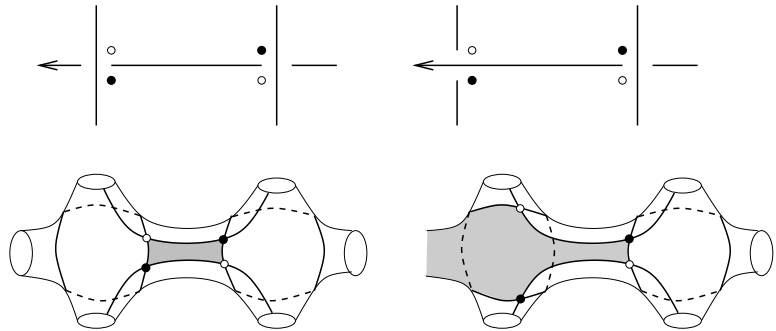}}
\caption{\label{fig:SpinCLemma}
{\bf{Illustration of Lemma~\ref{lemma:SpinCLemma}}}\qua The top row
represents the original projection diagram: the light circles represent
the state $x$, while the dark ones represent $y$. In the
second row, we have the corresponding Heegaard picture, with the
support of the homotopy class $\phi\in\pi_2(\x,\y)$ lightly shaded.
(Note that there are two other cases not pictured here, but the
corresponding pictures are the same as these, only viewed from
underneath.)}
\end{figure}

\begin{proof}
Orient the edge between $v_1$ and $v_2$ as it appears in the
knot. There are now four cases, according to whether $v_1$ or $v_2$
are under- or over-crossings.

Suppose that the edge takes us from under-crossing at $v_1$ to an
over-crossing at $v_2$ (i.e.\ our edge is on the bottom at $v_1$ and on the
top at $v_2$). In this case, we claim that after possibly switching
the roles of $\x$ and $\y$, we can find a homotopy class
$\phi\in\pi_2(\x,\y)$, with 
$n_{\BasePt}(\phi)=0$ and
$n_{\FiltPt}(\phi)=1$.  Indeed, the domain associated to $\phi$,
$\cald(\phi)$, has all local multiplicities zero or one;
topologically, it is a connected surface with a collection of circle boundary
components. Assume for a moment that there are no intermediate
crossings between $v_1$ and $v_2$. In this case, the topology of
$\cald(\phi)$ is given as follows: the genus of $\cald(\phi)$ is given
by the number of vertices encountered twice between $v_2$ and
the final  point, it has one boundary component corresponding to the meridian
$\mu$, it has one boundary circle for each vertex encountered only
once between $v_2$ and the final point (the corresponding
$\beta$-circle), and there is one final boundary component (with four
distinguished corner points) which is formed from arcs in
$\beta_{v_1}$, $\beta_{v_2}$, and the two $\alpha$-curves
corresponding to the two regions which contain the edge from $v_1$ to
$v_2$. This is pictured in Figure~\ref{fig:SpinCLemma}. In the case
where there are intermediate crossings between $v_1$ and $v_2$, the
surface looks much the same, except that now there are additional
circle components, one corresponding to each compact region in the
complement in $\R^2$ of the part of the projection between $v_1$ and
$v_2$ (or equivalently, one for each vertex between $v_1$ and
$v_2$). Specifically, these circle components are the $\alpha$-curves
of these intermediate regions.

The detailed description of the topology of $\cald(\phi)$ is not
particularly relevant to the proof of the present lemma (though it is
relevant in the proof of the next one); all we need here is the fact
that $n_\BasePt(\phi)=0$ and $n_\FiltPt(\phi)=1$, from which it
follows immediately that $\Filt(\x)-\Filt(\y)=1$. From the
definition of $S(x)$, it follows easily that $S(x)-S(y)=1$
(i.e.\ independent of the orientation on the two pieces of $K$
transverse to our edge at $v_1$ and $v_2$) as well, verifying the
lemma in this case.

Suppose that both $v_1$ and $v_2$ are under-crossings. When there are
no intermediate intersection points between $v_1$ and $v_2$, it is
easy to find a square which serves as $\cald(\phi)$ for some
$\phi\in \pi_2(\x,\y)$ with
$n_{\BasePt}(\phi)=n_{\FiltPt}(\phi)=0$, supported on the region of
the Heegaard surface corresponding to the edge, as pictured in
Figure~\ref{fig:SpinCLemma}.  In particular, the filtration
difference is zero. When there are intermediate intersection points,
the region $\cald(\phi)$ now is a square with a finite number of circles
removed, but it is still disjoint from $\BasePt$ and $\FiltPt$. Now,
it is easy to see that the formula $S(x)-S(y)$ gives zero, independent
of the crossing signs of the two vertices.

The remaining two cases (where $v_2$ is an over-crossing)
are handled similarly.
\end{proof}

\begin{lemma}
\label{lemma:GradingForConsecutives}
Suppose that the two states $x$ and $y$ differ by a transposition, at 
vertices $v_1$ and
$v_2$. Then, the absolute value of the
difference in gradings between the corresponding generators
$\x$ in $\y$ is one. More precisely, if
$\x$ is represented by light dots in
Figure~\ref{fig:GradingForConsecutives},
and $\y$ is represented by the dark ones, then $\gr(\x)-\gr(\y)=1$.

\begin{figure}[ht!]
\cl{\epsfxsize3.2in\epsfbox{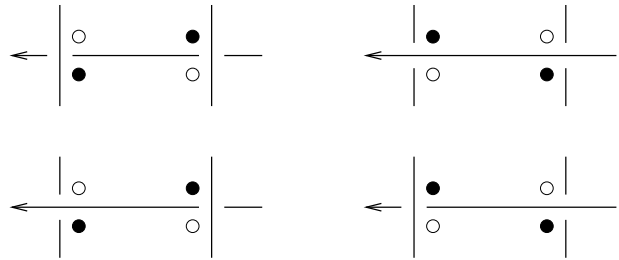}}
\caption{\label{fig:GradingForConsecutives}
{\bf{Grading difference for intersection points which differ by a transposition}}\qua 
Let $x$ be the state represented by the dark dots, and $y$
be the one represented by the light ones. Then, 
$\gr(x)-\gr(y)=1$.}
\end{figure}

\end{lemma}

\begin{figure}[ht!]
\relabelbox\small
\cl{\epsfxsize3in\epsfbox{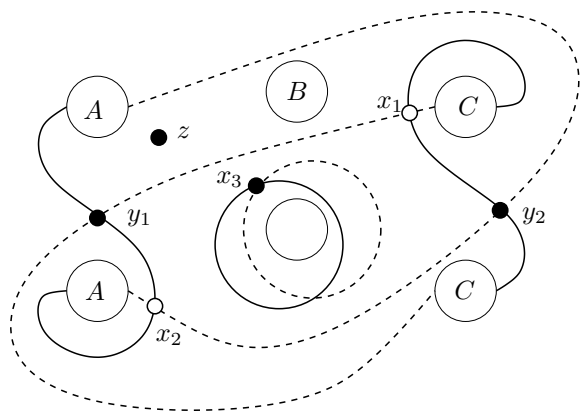}}
\adjustrelabel <-2pt, 0pt> {A1}{$A$}
\adjustrelabel <-2pt, 0pt> {A}{$A$}
\adjustrelabel <-2pt, 0pt> {B}{$B$}
\adjustrelabel <-2pt, 0pt> {C1}{$C$}
\adjustrelabel <-2pt, 0pt> {C}{$C$}
\adjustrelabel <-2pt, 0pt> {x1}{$x_1$}
\relabel {x2}{$x_2$}
\relabel {x3}{$x_3$}
\relabel {y1}{$y_2$}
\relabel {y2}{$y_1$}
\relabel {z}{$z$}
\endrelabelbox
\caption{\label{fig:MasCalc}
{\bf{Maslov index calculation of a square minus a disk}}\qua  Above, we
have pictured a $g=3$ Heegaard diagram for $\#^2(S^2\times S^1)$.  The
dark lines represent $\beta$-curves, and the dashed ones represent
$\alpha$-curves. Note that the disks labeled with capital letters
$A$, $B$, or $C$ are to be removed, and their boundaries are pairwise
identified, according to their pictured labeling. There are four
intersection points in $\Ta\cap\Tb$, all of which represent the
$\SpinC$ structure with trivial first Chern class, using the 
pictured reference point $z$.  Let
$\x=\{x_1,x_2,x_3\}$ and $\y=\{y_1,y_2,y_3\}$ with $y_3=x_3$.
It is easy to find a domain for a Whitney
disk $\phi\in\pi_2(\x,\y)$ which is a square with a disk
removed. Moreover, the grading difference between these two elements
is one (they are also connected by a square).}
\end{figure}

\begin{proof}
Domains with $n_\BasePt(\phi)=0$ connecting  intersection points which differ by a transposition
have already been demonstrated in the proof of
Lemma~\ref{lemma:SpinCLemma} (see Figure~\ref{fig:SpinCLemma}). Our
task here is to calculate their Maslov index (and, indeed, to show
that it is one in all the above cases). 

In all cases, we have seen that the domains $\cald(\phi)$ are
topologically a square with a finite number of handles attached, and a
finite number of disks removed. Letting $\x,\y\in\Ta\cap\Tb$ denote
the intersection points corresponding to the states $x$ and $y$, and
writing $\x=\{x_1,\ldots,x_{g}\}$, and $\y=\{y_1,\ldots,y_{g}\}$, we have
(for some numbering) that for $i=3,\ldots,g$, $x_i=y_i$, while $x_1$,
$y_1$, $x_2$, and $y_2$ are the four corner points of the square. Now,
in the interior of each handle attached to the square, we have a point
$x_i$ (for $i>2$), and also on each circle we have a point $x_i$, and all
the remaining intersection points $x_i$ (with $i>2$) lie in the exterior
of $\cald(\phi)$. Thus, these extra intersection points do not affect
the Maslov index. Moreover, each of the handles can be deleted, while
deleting the corresponding intersection point $x_i$ without affecting
the Maslov index. (This latter operation corresponds to destabilizing
the Heegaard diagram, see~\cite{HolDisk}.)

It remains then to calculate the Maslov index of a homotopy class
$\phi$ whose domain consists of square minus a finite collection of
disks. We claim that for such a homotopy class, $\Mas(\phi)=1$.  This
can be seen, for example, by finding such a homotopy class in a
suitably chosen (genus $n+2$) Heegaard diagram for $\#^{n+1}(S^1\times
S^1)$ (where here $n$ denotes the number of disks removed), and using
the fact that $\HFa(\#^{n+1}(S^1\times S^1),\spinc_0)\cong
H_*(T^{n+1};\Z)$ as relatively $\Z$-graded groups where here
$\spinc_0$ denotes the $\SpinC$ structure with first Chern
class equal to zero (cf.\ Section~\ref{HolDiskOne:sec:HandleSlides}
of~\cite{HolDisk}). Specifically, we can find such a domain which
connects two intersection points whose absolute gradings are known to
differ by one.  We illustrate this in the case where $n=1$, in
Figure~\ref{fig:MasCalc}.
\end{proof}

We find it useful to introduce one more notion before proceeding to
the proof of Theorem~\ref{thm:States}.  Recall that the decoration on
the knot projection (the distinguished edge, and the orientation on
the knot projection) induces a natural ordering on all the edges of
the knot projection. If $v$ is any vertex in the knot projection, the
{\em edge after $v$} is the last edge (in this ordering) whose
closure contains $v$. The edge which immediately precedes it will be called
the {\em $v$-penultimate edge}, and we denote it by $e_v$. Note also that the
orientation of $K$ and the reference point specifies an ordering on
the vertices, as follows. We say that $v_1<v_2$ if the vertex $v_1$ is
crossed the second time before $v_2$ is (on a path beginning 
in the interior of the distinguished edge and following $K$).

There is a {\em canonical state} $x_0$ which is
uniquely characterized by the property that for each crossing $v$,
$x_0(v)$ is one of the two quadrants whose closure contains $v$-penultimate
edge, according to the following lemma:

\begin{figure}[ht!]
\cl{\epsfxsize2in\epsfbox{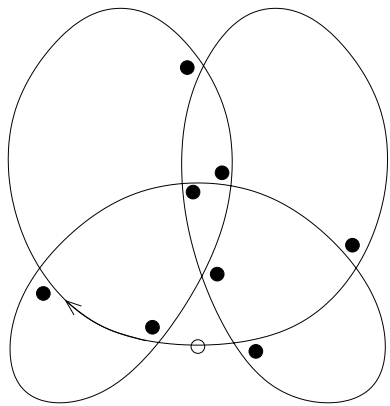}}
\caption{\label{fig:EightSixteen2}
{\bf{The canonical state $x_0$ for $8_{16}$}}\qua  We have illustrated
here the canonical state $x_0$ for the alternating knot $8_{16}$. The
light circle indicates the distinguished edge, and the arrow indicates
the orientation.}
\end{figure}

\begin{lemma}
\label{lemma:CanonicalState}
The canonical state $x_0$ is well-defined.
\end{lemma}

\begin{proof}
Let $R(G)$ denote the set of regions in $S^2-G-A-B$. Let $X_0=A\cup
B$. We will inductively define the canonical state $x_0$ by finding an
ordering $\{v_1,\ldots,v_m\}$ of all the vertices in the graph $G$, with
the property that $x_0|\{v_1,\ldots,v_n\}$ is uniquely defined (note that
it is {\em not} the same as the ordering induced by the orientation of
$K$). Correspondingly, we will exhaust $S^2$ by a sequence of regions
$$X_0\subset X_1 \subset \cdots\subset X_m=S^2$$ with the property that
$x_0|\{v_1,\ldots,v_n\}$ maps onto the set of regions in $X_n-G-X_0$ (and
hence $X_{n+1}$ is defined by attaching to $X_n$ a region $R_n\in
R(G)$).  In our induction hypothesis, we will assume that the
$v_i$-penultimate edge is contained in the interior of $X_i$.

For the initial step ($n=0$), $x_0$ is vacuously defined. 

For the inductive step, either $X_n=S^2$ (in which case we are
finished), or $X_n$ is a proper subset of $S^2$, in which case it must
have corners, since we have assumed that our graph belongs to a knot
projection. Consider the last corner point $v_{n+1}$ of $X_n$ (with
respect to the ordering induced by the orientation on the knot). It is
clear that the $v_{n+1}$-penultimate edge must appear in two quadrants, one
of which is contained in $X_n$, and the other of which is contained in a region
$R_{n+1}\subset S^2-X_n$ (for if this hypothesis is not satisfied, we
would simply be able to find a later corner vertex). It follows
then that $v_{n+1}\not\in
\{v_1,\ldots,v_n\}$ (for the edges before all those vertices are all contained
in the interiors of $X_n$). We then define $x_0(v_{n+1})$ to be the
quadrant of $R_{n+1}$ containing the $v_{n+1}$-penultimate edge. Let
$X_{n+1}=X_{n}\cup R_{n+1}$. It is now clear from the construction of
$X_{n+1}$ that the $v_{n+1}$-penultimate edge is contained in the interior
of $X_{n+1}$.

Observe that the above argument not only constructs the canonical
state $x_0$ but, since there was no ambiguity in the definition,
establishes its uniqueness.
\end{proof}

\proof[Proof of Theorem~\ref{thm:States}]
The correspondence between states and generators of $\CFKa(S^3,K)$ in
the appropriate Heegaard diagram was already explained in the
beginning of this section. We adopt here the notational convention
that if $x$ and $y$ are states, then $\x$ and $\y$ are their
corresponding intersection points.

To prove the assertion about filtration levels, we appeal to
Lemma~\ref{lemma:SpinCLemma}, according to which
if $x$ and $y$ are any two states which differ by a transposition,
then 
$$S(x)-S(y)=\Filt(\x)-\Filt(\y).$$ Combining this with Kauffman's
theorem, we see that there is a constant $c_1$ which {\em a priori}
depends on the knot projection, with the property that if $x$ is any
state, then
$$S(x)=\Filt(\x)+c_1.$$ In particular, it follows that for some suitable
choice of signs (given by the the $\Zmod{2}$
grading of the intersection point $\mathbf x$), the polynomial $$\Gamma_K (T)
=
\sum_{x\in {\mathfrak S}} (\pm 1) T^{S(x)},$$ has the form $\Gamma(T)=T^{c_1} \cdot
\Delta_K(T)$ where
$\Delta_K(T)$ denotes the symmetrized
Alexander polynomial. Indeed, the fact that that $c_1=0$ follows from
the fact that $\Gamma_K(T)$ is symmetric: indeed, it coincides with the
Conway-normalized Alexander polynomial, according to Chapter VI
of~\cite{KauffmanTwo}.

When $x$ and $y$ are states which 
differ by a transposition, and $\x$ and $\y$ their corresponding
intersection points, it follows easily from 
Lemma~\ref{lemma:GradingForConsecutives} that
$$\gr(\x)-\gr(\y)=M(x)-M(y).$$ From Kauffman's theorem, it then readily follows
that there is a constant $c_2$ with the property that
for any state $x$,
$$\gr({\bf x})=M(x)+c_2.$$ 
To see  that $c_2=0$, we 
verify that if $x_0$ is the canonical state
(whose existence was established in Lemma~\ref{lemma:CanonicalState} above)
and ${\bf {x_0}}$ its corresponding intersection point,
$$M(x_0)=0=\gr({\mathbf x_0}).$$ Indeed, it is straightforward to see that $M(x_0)=0$:
the local contributions $m(x_0,v)$ vanish for each vertex.

To see that ${\mathbf x_0}$ has vanishing absolute grading, we proceed
as follows.  One can reduce the Heegaard diagram for the knot
described above to another Heegaard diagram for $S^3$ by handlesliding
the $\beta$-curves of vertices.  We then arrange the $\beta_v$-curves
in descending order, according to this ordering of the corresponding
vertices, and then handleslide them ``forwards'' across the reference
point $\FiltPt$ (but never across $\BasePt$). In this manner, we
obtain a new Heegaard diagram for $S^3$, where the $\beta$-curve at
any vertex $v$ now meets only the up to two $\alpha$-curves
corresponding to the two quadrants which contain the $v$-penultimate
edge.  Thus, the canonical state also induces an intersection point
${\mathbf x}_0'$ for this new Heegaard diagram, and indeed its
absolute degree agrees with that of ${\mathbf x}_0$, since the
handleslides never crossed $\BasePt$ (compare Section~\ref{HolDiskOne:sec:HandleSlides} of~\cite{HolDisk}). Moreover, the uniqueness of
Lemma~\ref{lemma:CanonicalState} ensures that ${\mathbf x}_0'$
is, in fact, the only intersection point in
$\Ta\cap\Tb$ for this new Heegaard diagram, so it must have absolute
degree zero.
\qed

\section{Results on alternating knots}
\label{sec:AlternatingKnots}

For the purposes of Theorem~\ref{thm:KnotHomology}, it is useful to
have the following explicit description of the signature $\sigma(K)$
of an alternating link, which follows from work of Gordon and Litherland
\cite{GordonLitherland} as interpreted by 
Lee, see~\cite{EunSooLee}, which we now recall.

\begin{figure}[ht!]
\cl{\epsfxsize1in\epsfbox{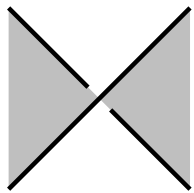}}
\caption{\label{fig:ColorConventions}
{\bf{Coloring conventions for alternating knots}}\qua 
We adopt the pictured convention in the statement of Theorem~\ref{thm:GoLi}.}
\end{figure}

Consider the four quadrants meeting at some vertex $v$.  If we orient
the boundary of one of these quadrants $Q$ (in the manner induced from
the orientation of the plane), we obtain an ordering on the two
consecutive edges contained in the boundary of $Q$. If the first of
these edges is part of an ``overcross'' at $v$, we call $Q$ an {\em
over-first} quadrant; otherwise, we call $Q$ an {\em under-first}
quadrant.  The alternating condition on a knot projection is
equivalent to the condition that all the over-first quadrants have the
same color (and hence that all the under-first ones have the other
color). For definiteness, we color all the under-first quadrants
white, as illustrated in
Figure~\ref{fig:ColorConventions}.

\begin{theorem}[Theorem~6 of~\cite{GordonLitherland}, 
see also Proposition~3.3 of~\cite{EunSooLee}]
\label{thm:GoLi}
Let $K$ be a knot with alternating projection $G$. Then, with the
coloring conventions illustrated in Figure~\ref{fig:ColorConventions},
the signature of $K$ is calculated by the formula
$$\sigma(K)=\#({\text{black regions in the planar projection}}) -
\#(\text{positive crossings})-1.$$
\end{theorem}

\proof[Proof of Theorem~\ref{thm:KnotHomology}]
Let $x$ be any state.
Glancing at the definitions of the local contributions $m(x,v)$ and $s(x,v)$,
one sees that
$$m(x,v)-s(x,v)+ \left(\frac{\epsilon(v)+1}{2}\right)
=
\left\{\begin{array}{ll}
0 &{\text{if $x(v)$ is an under-first quadrant}} \\
\OneHalf & {\text{if $x(v)$ is an over-first quadrant.}} 
\end{array}\right.
$$
Adding this up over all the vertices $v$, and
bearing in mind that the over-first quadrants are all black, and that 
there is exactly one black region (the distinguished one)
which is not represented as $x(v)$ for some vertex $v$, it follows that
\begin{gather*}
2(M(x)-S(x))=\hspace{3.8in}\\
\hspace{0.3in}\#({\text{black regions in the planar projection}}) -
\#(\text{positive crossings})-1,
\end{gather*}
a quantity which agrees with $\sigma(K)$, according to
Theorem~\ref{thm:GoLi}. The theorem now follows immediately from
Theorem~\ref{thm:States}.
\endproof

It is straightforward if tedious to verify that in fact the previous
theorem determines the filtered chain homotopy type of the knot
complex uniquely in terms of the Alexander polynomial and the
signature of the knot.  Rather than inflicting the necessary linear
algebra on our reader, we content ourselves here with the proof of
Theorem~\ref{thm:FloerHomology}.  For the proof, we will use the
relationship between the knot complex and the Floer homology of
three-manifolds obtained by sufficiently large surgeries on $K$, which
we recall presently.

Recall that when $K\subset S^3$ is a knot, we have the
$(\Z\oplus\Z)$-filtered complex $$\CFKinf(S^3,K)=\{[\x,i,j]\big|
\Filt(\x)+(i-j)=0\}.$$ Of course, the $\Z \oplus \Z$-grading of
$[\x,i,j]$ is $(i,j)$, and $$\gr[\x,i,j]=\gr(\x)+2i.$$ The subcomplex
where $i\leq 0$ (or $j\leq 0$) represents $\CFm(S^3)$, the whole
complex represents $\CFinf(S^3)$, and the quotient represents
$\CFp(S^3)$. More interestingly, this complex can be used to calculate
$\HFp(S^3_n(K))$ for sufficiently large $n$. Specifically, fix an
integer $s$, and let $\HFp(S^3_n(K),[s])$ denote $\HFp$ calculated
using a $\SpinC$ structure which extends over the two-handle
$B^4_n(K)$ to a $\SpinC$ structure satisfying $$\langle
c_1(\spinc),[{\widehat F}]\rangle = s-n,$$ where here $[{\widehat
F}]\in H_2(B^4_n(K);\Z)\cong \Z$ is a generator.  Consider the
quotient complex of $C$, which is generated by tuples $[\x,i,j]$ with
$i\geq 0$ or $j\geq s$, denoted $C\{i\geq 0~\text{or}~j\geq s\}$. 
It is shown in
Theorem~\ref{Knots:thm:LargePosSurgeries} of~\cite{Knots} that the homology 
of $C\{i\geq 0~\text{or}~j\geq s\}$ calculates $\HFp(S^3_n,[s])$,
under a map which shifts grading by $\frac{1-n}{4}$;
and indeed the natural projection map from $C\{i\geq 0~\text{or}~j\geq s\}$
to $C\{i\geq 0\}$ models the map induced from $B^4_n(K)$ induced by the 
$\SpinC$ structure $\spinc$. 

\proof[Proof of Theorem~\ref{thm:FloerHomology}]
The proof relies on the fact, which follows
from Theorem~\ref{thm:KnotHomology}, that
if we write $\sigma=\sigma(K)$, then for each intersection point $\x$
$$\gr(\x)=\Filt(\x)+\frac{\sigma}{2};$$ and hence, for each
$[\x,i,j]\in\CFKinf(S^3,K)$, 
\begin{equation}
\label{eq:GradingInTermsOfFilt}
\gr[\x,i,j]=i+j+\frac{\sigma}{2}.
\end{equation}
To calculate $\HFp(S^3_0(K),s)$, we calculate first $\HFp(S^3_n(K),s)$
for all sufficiently large $n$.

We have the following short exact sequence
$$
0\to
C\{\max(i,j-s)\geq 0\} 
\to
C\{i\geq 0\}\oplus C\{j\geq s\}
\to
C\{\min(i,j-s)\geq 0\}
\to0,
$$
where $H_*(C\{\min(i,j-s)\geq0\})$ is supported in degrees
$\geq s+\sigma/2$. 
It follows that
\begin{eqnarray}
H_{\leq s+\frac{\sigma}{2}-2}\left(C\{\max(i,j-s)\geq 0\}\right)
&\cong& H_{\leq s+\frac{\sigma}{2}-2}\left(C\{i\geq 0\}\oplus C\{j\geq s\}\right)
\nonumber \\
&\cong& \HFp_{\leq s+\frac{\sigma}{2}-2}(S^3)\oplus 
\HFp_{\leq -s+\frac{\sigma}{2}-2}(S^3) \nonumber \\
&=&\HFp_{\leq s+\frac{\sigma}{2}-2}(S^3) \label{eq:SmallDegrees}
\end{eqnarray}
(with the last equality following from the fact that $-s+\frac{\sigma}{2}-2<0$).

In the remaining degrees, we obtain information from  the 
short exact sequence
$$
0\longrightarrow
C\{\max(i,j-s)\leq -1\} \longrightarrow C \longrightarrow 
C\{\min(i,j-s)\geq 0\} \longrightarrow0
$$
Now, $H_*(C\{\max(i,j-s)\leq -1\})$ 
is supported in dimensions $\leq s+\sigma/2-2$.
Letting $R$ denote the the part in degree
$s+\sigma/2-2$, we have 
\begin{equation}
\label{eq:BigDegrees}
0\longrightarrow\HFinf_{\geq s+\frac{\sigma}{2}-1}(S^3) \longrightarrow H_{\geq s+\frac{\sigma}{2}-1}(C\{\min(i,j-s)\geq
0\})\buildrel\delta\over\longrightarrow R \longrightarrow \cdots
\end{equation}
Thus, there are two case for $\HFp(S^3_{n}(K))$ for large $n$.  When
$s+\frac{\sigma}{2}\leq 0$, the above arguments show that
$$H_*(C\{\max(i,j-s)\geq 0\})
\cong \HFinf_{\geq s+\frac{\sigma}{2}-1}(S^3)\oplus \Z_{(s+\frac{\sigma}{2}-1)}^{b_i},$$
for some non-negative integer $b_s$.
In fact, when $s+\frac{\sigma}{2}-1\geq 0$, the map from $C$ to
$C\{i\geq 0\}$ clearly factors through 
$C\{i\geq 0~{\text{or}}~j\geq s\}$. Thus, the above arguments show that
$$H_*(C\{\max(i,j-s)\geq 0\})
\cong \HFp(S^3)\oplus \Z^{b_s}_{(s+\frac{\sigma}{2}-1)}.$$
Ths identification of $H_*(C\{\max(i,j-s)\})\cong \HFp(S^3_n,[s])$ for
sufficiently large $n$, together with the integer surgery long exact
sequence (\cite{HolDiskTwo}, see also~\cite{AbsGraded}) $$
\cdots \longrightarrow\HFp(S^3) \buildrel{F_1}\over\longrightarrow \HFp(S^3_0(K),s) \buildrel{F_2}\over\longrightarrow\HFp(S^3_n(K),[s])
\buildrel{F_3}\over\longrightarrow\cdots
$$ 
now gives the result when $s\neq 0$.

In the case where $s=0$, we still have
a (possibly trivial) cyclic summand $\Z^{b_0}$ in
$\HFp(S^3_n(K),0)$ (for $n$ sufficiently large) supported in dimension
$\frac{1-n}{4}+\frac{\sigma}{2}-1$, with the property that
$$\HFp(S^3_n(K),0)\cong
\InjMod{-2\lceil\frac{-\sigma}{4} \rceil+\frac{1-n}{4}}\oplus \Z^{b_0}.$$
Now, in the integral surgeries long exact sequence,
the map $F_3$ consists of a sum of terms, 
each of which decreases the absolute grading 
by at least $\frac{1-n}{4}$. It follows immediately (again, using our
hypothesis that $\sigma<0$) that this map in the present case must
vanish.

The map $F_2$ now shifts degree by $(n-3)/4$ and the map $F_1$ shifts
degree by $-1/2$ (cf.\ Lemma~\ref{AbsGraded:lemma:CalcDegrees}
of~\cite{AbsGraded}), so the calculation of $\HFp(S^3_0(K),0)$ follows.

The only remaining piece now is the verification of
Equation~\eqref{eq:GetBs}. But this follows immediately from
Theorem~\ref{HolDiskTwo:thm:EulerOne} of~\cite{HolDiskTwo}, where the
Euler characteristic of $\HFp(Y_0)$ is identified with the torsion of
$S^3_0$, or, more precisely, provided that $i\neq 0$,
$$-t_i(K)=\chi(\HFp(S^3_0(K),i))$$
(see also~\cite{HolDiskTwo} for the statement when $i=0$).
\endproof

We now turn to the proofs of the corollaries listed in the introduction.

\proof[Proof of Corollary~\ref{cor:AlternateEstimate}]
This is an immediate application of the theorem (after reflecting $K$ if necessary), 
bearing in mind that, of course, $b_s\geq 0$.
\endproof

\proof[Proof of Corollary~\ref{cor:Froyshov}]
In general, we have that $d(S^3_1(K))=d_{\OneHalf}(S^3_0(K))-\OneHalf$
(see Proposition~\ref{AbsGraded:prop:CorrTermEquality}
of~\cite{AbsGraded}).  In the case where $\sigma(K)\leq 0$, the result
then is an immediate application of
Theorem~\ref{thm:FloerHomology}. When $\sigma(K)>0$, let $r(K)$ denote
the reflection of $K$; then we have (see~\cite{AbsGraded}) that
$$d_{\OneHalf}(S^3_0(K))=-d_{-\OneHalf}(S^3_0(r(K)))=\OneHalf$$ (with
the last equation following once again from
Theorem~\ref{thm:FloerHomology}).
\endproof

For the proof of Corollary~\ref{cor:AltContact}, we use the results
from ~\cite{HolDiskContact}, which in turn rely on results of
Giroux~\cite{Giroux}. Specifically, if $K\subset Y$ is a fibered knot
of genus $g$, then $\HFKa(-Y,K,-g)\cong \Z$. The image of the
generator of this group inside $\HFa(-Y)$ is shown
in~\cite{HolDiskContact} to depend on the knot only through its
induced contact structure $\xi_K$, giving rise to an element
$c(\xi_K)\in \HFa(-Y)$.  Moreover, when $Y\cong S^3$ (or, more
generally, $\xi$ is any contact structure over a three-manifold $Y$
whose whose induced $\SpinC$ structure has torsion first Chern class),
$c(\xi)$ is a homogeneous element whose absolute degree $c(\xi)$
coincides with the Hopf invariant of $\xi$. Finally, in
Theorem~\ref{HolDiskContact:intro:OverTwisted}
of~\cite{HolDiskContact}, the invariant is shown to vanishing for
overtwisted contact structures.

\proof[Proof of Corollary~\ref{cor:AltContact}]
According to Theorem~\ref{thm:KnotHomology} the degree of an element
in filtration degree $-g$ (and hence, as above, the Hopf invariant of
the induced homotopy class of two-plane field) is given by
Equation~\eqref{eq:EquationForHopf}. Note that the sign appearing in 
front of the signature occurs because, in the definition of $c(\xi)$,
we are reverse the orientation on the ambient three-manifold, which
is equivalent to reflecting the knot.

In the case where this Hopf invariant vanishes, the induced element in
$\HFa(S^3)$ must be non-trivial, for it is the only generator in
degree zero (again, according to Theorem~\ref{thm:KnotHomology}).
Thus, according to Theorem~\ref{HolDiskContact:intro:OverTwisted}
of~\cite{HolDiskContact}, the induced contact structure is tight.
\endproof

\proof[Proof of Corollary~\ref{cor:AltContact2}]
The fact that the Hopf invariant of is non-positive, follows readily
from Equation~\eqref{eq:EquationForHopf}, together with the fact that
$|\sigma(K)|\leq 2g$.

Moreover, we have seen in Corollary~\ref{cor:AltContact} that the overtwisted
contact structure with vanishing Hopf invariant cannot be represented
by an alternating knot; while it is clear that unknot represents the
tight contact structure.

Finally, using Equation~\eqref{eq:EquationForHopf}, we see that
for each $i>0$,
the $i$-fold connected sum of the figure eight knot realizes the
contact structure with Hopf invariant $-i$.
\qed

\section{Alternating links}
\label{sec:AltLinks}

Recall that in~\cite{Knots}, we defined a generalization of the knot
invariants $\HFKa$ to the case of links. These link invariants satisfy
a skein exact sequence
(cf.\ Theorem~\ref{Knots:thm:SkeinExactSequence} of~\cite{Knots}):
suppose that $p$ is a positive crossing for a projection of a link $L_+$,
for which both strands belong to the same component of $L_+$, then there is a
long exact sequence (for each $s\in\Z$) of the form: 
\begin{equation}
\label{eq:SkeinExactSequence}
\cdots\longrightarrow
\HFKa(L_-,s)\buildrel{f}\over\longrightarrow
\HFKa(L_0,s)\buildrel{g}\over\longrightarrow
\HFKa(L_+,s)\longrightarrow\cdots
\end{equation}
where $L_-$ is the modified version of $L_+$ (with a crossing-change
at $p$), and $L_0$ is the resolution at $p$ of $L_+$. Both maps $f$
and $g$ drop absolute grading by $1/2$, where the remaining 
map is non-increasing on the absolute grading.

For the following statement, recall that a link $L$ called a {\em
non-split, alternating link} if it has a projection which is
connected, and also, if we traverse any component of $L$, the
crossings in this projection alternate between over-crossings and under-crossings.

\begin{theorem}
\label{thm:LinkHomology}
Let $L\subset S^3$ be a non-split, oriented, alternating
link  in the three-sphere,
and let $\Delta_L$ be its Alexander-Conway polynomial. Writing
$$(T^{-1/2}-T^{1/2})^{n-1}\cm \Delta_L = a_0+\sum_{s>0} a_s
(T^s+T^{-s}),$$ we have that $\HFKa(S^3,L,s)$ is supported entirely in
dimension $s+\frac{\sigma}{2}$, and indeed $$\HFKa(S^3,L,s)\cong
\Z^{|a_s|}.$$ Here, $\sigma$ is the signature of the link $L$.
\end{theorem}

\begin{figure}[ht!]
\cl{\epsfxsize4in\epsfbox{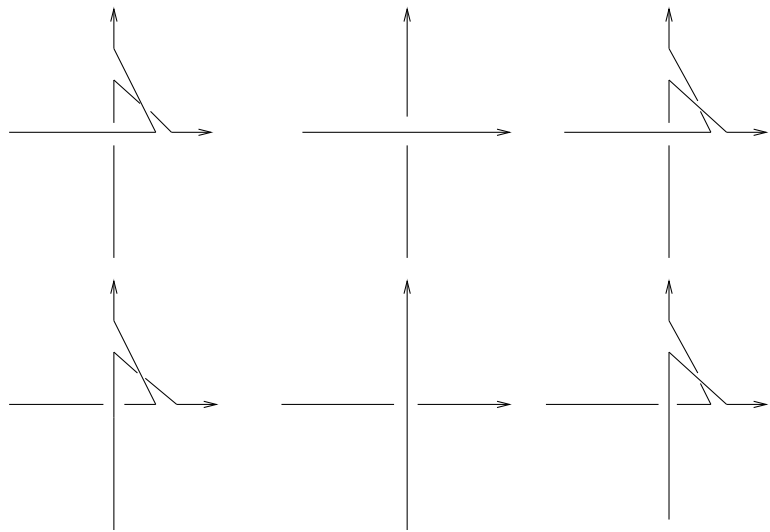}}
\caption{\label{fig:AlternateChange}
{\bf{Skein moves on alternating links}}\qua  On the left, we have two
possible candidates for $L_-$, in the middle we have the two versions
of $L$, while on the right we have two candidates for $L_+$. It is
easy to see that if the links $L$ are alternating, then the changes
$L_-$ and $L_+$ can also be arranged to alternate (after cancelling an
extra pair of intersection points, if necessary.}
\end{figure}

\begin{proof}
Recall first that the skein exact sequence 
can be used to show that 
$$\chi(\HFKa(S^3,L,i))=a_i$$
(cf.~\cite{Knots}).

In view of this calculation, the theorem is obtained by induction on
the number of components of $L$, with Theorem~\ref{thm:KnotHomology}
as base case.

For the inductive step, let $p$ be an intersection point where two
different strands of $L$ meet. We can find two links links $L_-$ and
$L_+$ with one more intersection point $q$, both of which admit
alternating projections, and which have the the property that their
resolution $L_0$ at $q$ is our original $L$. The two cases, according
to the sign of the intersection point $p$, are illustrated in 
Figure~\ref{fig:AlternateChange}. 

When $p$ is a positive intersection point for $L$, we see
that (after the obvious cancellation), $L_-$ has one fewer positive
intersection points than $L$ does, while $L_+$ has one more positive
intersection point.  Moreover, the number of black regions (using the
coloring conventions of Figure~\ref{fig:ColorConventions}) are the
same for all three links. If, on the other hand, $p$ is a negative
intersection point, then the number of black regions for $L_-$ is one
greater than number for $L$, which in turn is one greater than the
number for $L_+$. Moreover, the number of positive intersection points
is the same for all three. Thus, applying Theorem~\ref{thm:GoLi}, we
can conclude that in either case,
\begin{equation}
\label{eq:SignatureChange}
\sigma(L_-)-1=\sigma(L)=\sigma(L_+)+1.
\end{equation}
It is now straightforward to conclude the result for $L=L_0$ from the
skein exact sequence, and the inductive hypothesis on $L_-$ and $L_+$.
\end{proof}

\begin{figure}[ht!]
\cl{\epsfxsize4.5in\epsfbox{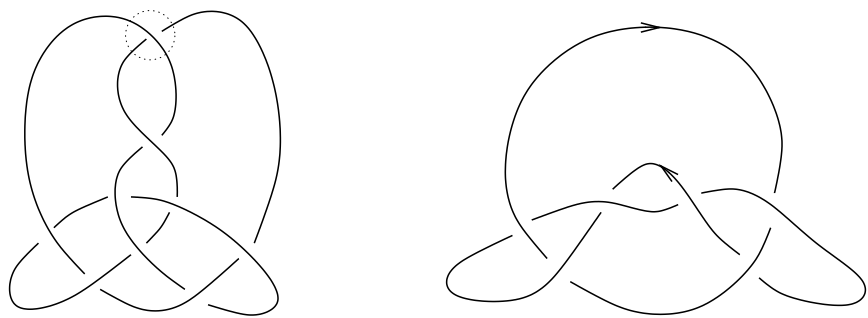}}
\caption{\label{fig:Nine48}
{\bf{The knot $9_{48}$}}\qua We have illustrated this nine-crossing
knot. If the crossing circled with a dotted circle is switched, we
obtain the right-handed trefoil; while if the crossing is resolved, it
is easy to see that we obtain the (oriented) alternating link pictured
on the right.}
\end{figure}

Theorem~\ref{thm:LinkHomology} can be used to give easy
generalizations to (non-split) alternating links of the results stated
in the introduction for alternating knots.  Rather than enumerating
these, we use Theorem~\ref{thm:LinkHomology} to give a calculation of
$\HFKa$ for the (non-alternating) knot pictured on the left in
Figure~\ref{fig:Nine48} (known as ``$9_{48}$'' according to the
standard knot tables, cf.~\cite{BurdeZieschang}).

If we change the indicated crossing, we obtain the
right-handed trefoil
$K_-$, which has 
\begin{eqnarray*}
\sigma(K_-)&=&-2 \\
\Delta_{K_-}&=& T^{-1}-1+T.
\end{eqnarray*}  If the knot crossing is
resolved, we obtain the two-component link $L$ pictured on the right in
Figure~\ref{fig:Nine48} (given the specified orientation).
It is straightforward to calculate that 
\begin{eqnarray*}
\sigma(L)&=& -1 \\
(T^{-1/2}-T^{1/2})\cm \Delta_L &=& T^{-2}-6T^{-1}+10-6T+T^2.
\end{eqnarray*}
It is now an immediate application of the skein exact sequence
and Theorem~\ref{thm:LinkHomology}
that the conclusion of Theorem~\ref{thm:LinkHomology} holds for $9_{48}$
(and hence also the conclusion of Theorem~\ref{thm:FloerHomology}),
even though $9_{48}$ does not possess an alternating projection.

More calculations of knot homology groups are given~\cite{calcKT}.

\end{document}